\documentclass[12pt]{article}
\begin{document}
\begin{center}
\large{\bf{Uniform distribution on the sphere and caps}}
\bigskip

{\large{\bf{Aljo\v{s}a Vol\v{c}i\v{c}}}}

Universit\`a della Calabria
\end{center}
\bigskip
\bigskip

\begin{abstract}
In this note we will consider the question when from the appropriate behavior  of a sequence of points on caps we can conclude that the sequence is uniformly distributed on the sphere.
\end{abstract}

\section {\bf Introduction} 

We denote by $S^{n-1}$ the unit sphere, the set of all unit vectors of the euclidean space $I\!\!R^n$. By $u\cdot v$ we denote the usual scalar product in $I\!\!R^n$. A {\it cap} is a portion of the sphere cut off by a hyperplane.

Suppose $Q$ is a regular Borel probability on $S^{n-1}$. We say that a sequence of unit vectors $\{u_j\}$ is $Q$-uniformly distributed on $S^{n-1}$, if 
$$\lim_{N\rightarrow \infty}\frac{1}{N}\sum_{j=1}^N \chi_A(u_j)=Q(A)  \eqno (1)$$
holds for all  $Q$-continuity sets, i.e. sets $A$ such that $Q(\partial A)=0$, where $\partial A$ is the relative boundary of $A\subset S^{n-1}$.
\smallskip

When $P$ is the probability corresponding to the normalized Hausdorff measure ${\cal H}^{n-1}$ on the sphere, we usually simply speak about {\it uniform distribution}, without mentioning the probability $P$.
\bigskip

Peter Gruber posed during the conference on Convex and Discrete Geometry (July  2009, Vienna)  a problem which can be reformulated in the following way:
\bigskip

\noindent {\bf Problem 1} Suppose $C\subset I\!\!R^n$ is a cap having a fixed size, $P(C)=a$, say, with $0<a<\frac{1}{2}$. Is it true that if
$$\lim_{N \rightarrow \infty} \frac{1}{N}\sum_{k=1}^N \chi_{\rho C}(u_j)=a \eqno (1)$$
holds for every rotation $\rho$ of the cap $C$, then the sequence $\{u_j\}$ is uniformly distributed on $S^{n-1}$?
\smallskip

Using the terminology of the theory of uniform distribution, this ammounts to ask whether the family of all caps of given size is a {\it discrepancy system} ([DT], [KN]).
\bigskip

We shall first give a simple (and complete) solution of the problem for $n=2$ and then provide an answer for the case $n\ge 3$.

\section {\bf The planar case} 

We will solve in this section the planar version of Problem 1. A cap here is an arc and the probability $P$ is just the normalized arc-lenght.

The following theorem shows that the answer depends on the size of the cap. 
\bigskip

\noindent {\bf Theorem 1} {\it If (1) holds for every rotation $\rho$ of the arc $C$, then

(i) if  $a$ is rational, $\{u_j\}$ need not be uniformly distributed;

(ii) if $a$ is irrational, $\{u_j\}$ has to be uniformly distributed.}
\smallskip

\noindent {\bf Proof} 

(i) Suppose $a=\frac{r}{s}$ with $r$ and $s$ coprime positive integers.  Put $F(\vartheta)=\frac{1}{2}\sin2s \vartheta$. Since $F$ is odd and $1+F$ is positive,
$$Q(B)=\frac{1}{2\pi}\int_B (1+F(\vartheta))\, d\vartheta$$
is a probability on the Borel subsets of $S^1$, distinct from $P$. By a generalization of von Neumann's theorem ([vN], [N]), there exists a $Q$-uniformly distributed sequence $\{v_j\}$ on $S^1$. Since the integral of $F$ on every arc of length $2\pi a$ is zero, equation (1) holds. But if $A$ is a $P$ (and $Q$) -continuity set for which $P(A)\not =Q(A)$, we have
$$\lim_{N \rightarrow \infty} \frac{1}{N}\sum_{j=1}^N \chi_A(v_j)=Q(A)\,,$$
showing that $\{v_j\}$ is not $P$-uniformly distributed.
\bigskip

(ii)) Suppose now that $a$ is irrational. Let us recall that a sufficient condition for uniform distribution of $\{u_j\}$ is that 
$$\lim_{N \rightarrow \infty} \frac{1}{N}\sum_{j=1}^N \chi_{[0, \beta]}(u_j)=\frac{\beta}{2\pi}\eqno (2)$$
holds for every $\beta$ belonging to a dense subset of $[0, 2\pi[$.

Since $a$ is irrational, by a well known result due to Kronecker (improved by Weyl), the set $\{2m\pi a: m\in I\!\!N\}$ is dense on $S^1$.

Let then $2m\pi a =\beta+2k\pi$, with $0<\beta <2\pi$ and let us evaluate
$$\frac{1}{N}\sum_{j=1}^N\sum_{i=1}^m \chi_{[(i-1)2\pi a, i2\pi a]}(u_j)=k+\frac{1}{N}\sum_{j=1}^N\chi_{[0,\beta]}(u_j)\,.$$

Note now that the left-hand side tends, by (1), to $ma=\frac{\beta}{2\pi}+k$. Therefore (2) holds
and this concludes the proof.

\section {\bf  Higher dimension} 

In this section we will show that when $n\ge3$, condition (1) does not imply, in general, that the sequence $\{u_j\}$ is uniformly distributed. The idea of the proof is the same we used in proving part (i) of Theorem 1, but in higher dimension ``bare hands" are not enough and we have to use Theorem 2 (called the {\it freak theorem}) due to Ungar [U] for $n=3$ and generalized to higher dimension by Schneider in [S1] and [S2].
\bigskip

We need some notation. Suppose $F$ is a continuous odd function on $S^{n-1}$. For any $s\in ]0,1[$ let us define the function $\tau_s(x)=1$ if $x\ge s$ and $\tau_s(x)=0$ otherwise. Consider now the integral transform
$${\cal T}_s(F)(u)=\int_{S^{d-1}}\tau_s(u\cdot v)F(v)\, d P(v)\,.$$
\smallskip

If we denote by $C_s(u)$ the cap $\{v: u\cdot v \ge s \}$, the previous definition can be rewritten as 
$${\cal T}_s(F)(u)=\int_{C_s(u)}F(v)\, d P(v)\,.$$
\smallskip

\noindent {\bf Theorem 2} {\it Let $F$ be a continuous odd function on $S^{n-1}$. The condition ${\cal T}_s(F)(u)\equiv 0$ implies that $F\equiv 0$ if and only if $s$ is not a zero of a $d+2$ dimensional Legendre polynomial of even degree.} 
\smallskip

A closer look to the previous theorem (see for instance [G], Remarks and References to 3.5) tells us that the countable set $D$ of the zeros $s$ of such polynomials  is dense in $]0,1[$.

Note also that since $F$ is odd, 
$$\int_{S^{d-1}} F(v)\, d P(v)=0\,.\eqno (4)$$
\smallskip

We are now in position to prove the following result.
\bigskip

\noindent {\bf Theorem 3 } {\it If $s \in D$ then $\{u_j\}$ need not be uniformly distributed even if
$$\lim_{N \rightarrow \infty} \frac{1}{N}\sum_{n=1}^N \chi_{C_s(u)}(u_j)=P(C) \eqno (5)$$
holds for any spherical cap $C_s(u)$.} 
\smallskip

\noindent {\bf Proof}  
Let $s \in D$ and let $F$ be an odd continuous function (not identically zero) on $S^{d-1}$ satisfying (4).

 Since $F$ is bounded, for an appropriate coefficient $c$, $1+c\, F>0$ on the sphere and therefore
$$Q(B)=\int_B (1+c\, F(u)) \,dP(u)$$
is a probability on the class of all Borel subsets on $S^{n-1}$, which is absolutely continuous with respect to $P$.  Obviously $P=Q$ on all the caps $C_s(u)$, but there are $P$ (and hence $Q$) -continuity sets $A$ such that $P(A)\not=Q(A)$.

By von Neumann's theorem recalled in Theorem 1, there exists a sequence $\{v_j\}$ of unit vectors which is uniformly distributed with respect to $Q$. This implies that 
$$\lim_{N \rightarrow \infty} \frac{1}{N}\sum_{j=1}^N \chi_A(v_j)=Q(A)\not =P(A))\,,$$
showing that $\{v_j\}$ is not $P$-uniformly distributed eventhough (5) holds, as we wanted to prove.
\bigskip

It would be of course interesting to complement this statement, proving also in dimension higher than two an analogue of part (ii) of Theorem 1, or clearifying anyway what happens if $s \not \in D$.
\bigskip

We can formulate this open question introducing a concept which has been studied in measure theory and probability, namely the uniqueness family of measurable sets.

A family of continuity sets $\cal U$ is said to be a {\it uniqueness} family for the probability $P$ if the knowledge of the values of $P$ on $\cal U$ determines $P$ uniquely.

The previous theorem  can be reformulated saying that a discrepancy system (introduced in Section 1) is a uniqueness family.
\bigskip

\noindent {\bf Problem 2} Is a uniqueness family for the normalized $n-1$ Hausdorff measure on $S^{n-1}$ allways a determining family? 
\bigskip

In measure theory there are many issues of the problem on wether the knowledge of the vaule of the measure on certain sets determines the measure uniquely (see for instance [D], [H], [R], {PT}, [Ro]), but this one does not seem to be classifiable among them.
\bigskip

 \noindent {\bf Acknowledgment} 

The author wishes to thank Professor Dietrich K\"olzow for helpful discussions concerning in particular the freak theorems.
\bigskip\bigskip\bigskip

\noindent {\bf References}
\bigskip\bigskip

\noindent [D] Davies R. O., Measures not approximable or not specifiable by means of balls. Mathematika, {\bf 18} (1971), 157-160 \bigskip

\noindent [DT], Drmota, M., Tichy R. F., {\it Sequences, discrepancies and applications.} Lecture Notes in Math. {\bf 1651} Springer, 1997
\bigskip

\noindent [G], Groemer H., {\it Geometric application of Fourier series and spherical harmonics.} Cambridge University Press, 1996
\bigskip

\noindent [H] Hoffmann-J¿rgensen J., Measures which agree on balls. Math. Scand. {37, 319Ð326 (1975) \bigskip

\noindent [KN] Kuipers L., Niederreiter H., {\it Uniform distribution of sequences.} Wiley and Sons, 1974
\bigskip

\noindent [N] Niederreirer H., A general rearrangement theorem for sequences. {\it Arch. Math.} {\bf 43} (1984), 530-534
\bigskip

\noindent [PT] Preiss D., Ti\v{s}er L., Measures in Banach spaces are determined by their values on balls, Mathematika, {\bf 38} (1991), 391-397
\bigskip

\noindent [R] Rana I. K., Determination of probability measures through group actions. {\it Z. Warscheilichkeitstheorie verw. Gebiete} {\bf 53} (1980), 197-206
\bigskip

\noindent [Ro]  Roberts M. L., Real Anal. Exchange {\bf 28}, no. 2 (2002), 635-640
\bigskip

\noindent [S1] Schneider R., Functions on a sphere with vanishing integrals over certain subspheres. {\it Journal of Mathematical Analysis and Applications} {\bf 26} (1969), 381-384
\bigskip

\noindent [S2] Schneider R.,  \"Uber eine Integralgleichung in der Theorie der konvexen K\"orper. {\it Math. Nachr.} {\bf 44} (1970), 55-75
\bigskip

\noindent [U] Ungar P., Freak theorem about functions on a sphere. {\it J. London Math. Soc.} {\bf 29} (1954), 100-103
\bigskip

\noindent [vN] von Neumann J., Gleichm\"assig dichte Zahlenfolgen. {\it Mat.
Fiz. Lapok} {\bf 32} (1925), 32-40 

\end{document}